
\documentstyle{amsppt}
\define\PP{{\Bbb P}}
\define\LL{\Cal L}
\predefine\Odanish\O
\redefine\O{\Cal O}
\define\E{\Cal E}

\define\A{\Cal A}
\define\m{\frak m}
\define\ord{\operatorname{ord}}
\define\codim{\operatorname{codim}}

\define\wt{\operatorname{wt}}
\define\Span{\operatorname{Span}\nolimits}
\define\Sym{\operatorname{Sym}}

\define\Proj{\operatorname{Proj}}
\define\im{\operatorname{Im}}
\define\Op#1{{\Cal O}_{\Bbb P^{#1}}}
\define\supp{\operatorname{supp}}
\define\lspan#1{\left<#1\right>}
\define\idsheaf#1{\Cal I_{#1}}
\define\set#1{\left\{\,#1\mid}
\define\endset{\,\right\}}
\define\Plucker{Pl\"u\discretionary{k-}{k}{ck}er}
\define\quasi{quasihomogeneous}
\define\emph#1{{\it #1\/}}

\tolerance 9999
\emergencystretch 3em
\hfuzz .5pt

\document

\topmatter
\author
       S.L'vovsky
\endauthor
\title On number of hypersurfaces containing projective
       curves
\endtitle
\date
      March 31, 1995
\enddate
\address
Serge L'vovsky,
Profsoyuznaya ul.,
20/9 kv.162.
Russia, 117292, Moscow.
\endaddress

\email
serge\@ium.ac.msk.su
\endemail
\abstract
Generalizing a classical lemma of Castelnuovo, we characterize
rational normal curves (resp.\ linearly normal elliptic curves)
as curves $C\subset \PP^n$ such that the number of linearly
independent hypersurfaces $Z\supset C$ of given degree~$m$
is maximal (resp.\ next to maximal).
\endabstract

\thanks
This research was partially supported by ISF grant MSC000
and Russian Foundation for Fundamental Research grant
95--01--00364
\endthanks

\endtopmatter

\document

\head 1. Introduction
\endhead
A well-known lemma of Castelnuovo (see, for
instance,~\cite{GH, Ch.~4, Sec.~3}) says that if $n(n-1)/2$
linearly independent quadrics pass through $d\ge 2n+3$ points in
general position in $\PP^n$, then these points lie on a rational
normal curve. As an immediate consequence, one can see that
\proclaim{Theorem 1.1 (Castelnuovo)}
If $n(n-1)/2$ linearly independent quadrics pass through a
non-degenerate curve $C\subset \PP^n$, then $C$ is a rational
normal curve.
\endproclaim
By induction on dimension, one can derive the following result
from this theorem:

\proclaim{Theorem 1.2}
If $X\subset \PP^n$ is a non-degenerate irreducible projective
variety and $c=\codim X$, then $X$ is contained in at most
$c(c+1)/2$ linearly independent quadrics; the equality is
attained iff the $\Delta$-genus of $X$ is~$0$.
\endproclaim
The condition ``$\Delta$-genus of $X$ is zero'' means that $\deg
X=c+1$; see \cite{H, page 51} for complete classification of
such varieties.

In an unpublished paper~\cite{Z}, F.L.~Zak proposed a new proof
and a generalization of this fact. The aim of this paper is to
give still another proof and a generalization of some results
of~\cite{Z}. The methods of proof are entirely different.

Theorem~1.2 can be generalized to the case of
$\Delta$-genus~$1$ (cf.~\cite{Z}):
\proclaim{Theorem 1.3}
Suppose that $X\subset \PP^n$ is a non-degenerate and
non-singular in codimension~$1$ irreducible projective variety
and put $c=\codim X$. If $X$ is contained in $c(c+1)/2-1$
linearly independent quadrics, then $\deg X=\codim X$ and smooth
one-dimensional linear sections of $X$ are linearly normal
elliptic curves.
\endproclaim
The varieties whose one-dimensional linear sections are linearly
normal elliptic curves are also completely classified (see, for
instance,~\cite{F}).

The main step in the proof of these theorems is to establish them
for curves (the rest is done by easy induction). If we consider
hypersurfaces of degree $>2$, the induction step does not work,
but some results for curves still can be obtained (see
section~1.1 for notation and terminology):
\proclaim{Theorem 1.4}
If $X\subset \PP^n$ is a non-degenerate irreducible projective
curve and $m \ge2$, then $h^0(\idsheaf X(m))\le \binom
{m+n}n-mn-1$. The equality is attained iff $X$ is a rational
normal curve.
\endproclaim
\proclaim{Theorem 1.5}
Suppose that $X\subset \PP^n$, $n\ge 2$, is a smooth
non-degenerate irreducible projective curve and $m \ge2$. If
$h^0(\idsheaf X(m)) < \binom {m+n}n-mn-1$, then $h^0(\idsheaf
X(m)) \le \binom{m+n}n - m(n+1)$. The equality is attained iff
$\deg X=n+1$ and the genus of $X$ is~$1$.
\endproclaim

The proof is based on the study of inflection points of the curve
$X$. It involves two ideas. The first (and very simple) one is
that the presence of an inflection point imposes an upper bound
on the number of hypersurfaces of degree $m$ containing $X$
(Proposition~3.1). The second idea is that that if
this bound is attained, then, under some additional conditions,
$X$ can be deformed to a certain monomial rational curve.
Actually, we prove that Hilbert polynomials of $X$ and the
monomial curve are the same rather than construct an explicit
deformation (Proposition~4.11).

Homogeneous rings of monomial curves, which we use in this paper,
have been studied by many authors (cf.\ \cite{B}, \cite{Ho}, and
references therein). Our Proposition~4.1 is contained in
Section~4 of \cite{B}. A computation of Hilbert polynomials of
monomial curves is contained in~\cite{B}, too: our
Proposition~4.9 agrees with Theorem~4.1 of~\cite{B}, while
methods of proof differ.

\subhead  Acknowledgements
\endsubhead
I am greatly indebted to F.L.Zak for his concern and numerous
helpful suggestions and discussions. I am grateful to S.~Lando
for useful discussions. This work was inspired by a reading of
the preprint~\cite{L}; I would like to thank J.M.~Landsberg for
an interesting email discussion of this preprint.

\subhead 1.1 Notation and conventions
\endsubhead
We work over an algebraically closed field of arbitrary
characteristic.

If $a$ and $b$ are integers, then $(a;b)=\set{n\in \Bbb
Z} a< n <b\endset$ and $[a;b]=\set{n\in \Bbb Z} a \le
n \le b\endset$.

For any subset $X\subset \PP^n$, denote by $\lspan X$ its linear
span. A projective subvariety $X\subset \PP^n$ is called
\emph{non-degenerate} iff $\lspan X=\PP^n$.

By $v_m\colon \PP^n \to \PP^{\binom{m+n}{n}-1}$ we denote the
\hbox{$m$-th} Veronese mapping.

If $0<a_1<\dots<a_n$ is an increasing sequence of positive
integers, then $\left<a_1\dots a_n \right>\subset \Bbb Z$
denotes the semigroup generated by $0, a_1,\dots, a_n$. If a
positive integer is not contained in this semigroup, it is called
its \emph{gap}.

\head 2. Preliminaries on inflection points
\endhead
In this section we chiefly recall some well-known definitions
(cf.~\cite{EH}). Let $C$ be a smooth projective curve, $\LL$ a
line bundle on $C$, $V\subset H^0(\LL)$ a linear subspace. We
will refer to the pair $(\LL, V)$ as a \emph{linear system} on
$C$.

Put $\dim V=n+1$. For any point $p\in C$, there exists a basis
$s_0,\dots,s_n$ of $V$ such that $\ord_p(s_0)< \ord_p(s_1)< \dots
<\ord_p(s_n)$. Denote $\ord_p(s_i)$ by $a_i(p)$. The sequence
$a_0(p), a_1(p), \dots, a_n(p)$ is called \emph{the vanishing
sequence at $p$}. For any $s\in V$, the number $\ord_p(s)$
coincides with one of the $a_j(p)$'s. The number $w(p)=\sum_{0\le
i\le n}(a_i(p)-i)$ is called \emph{weight} of $p$ with respect to
the linear system $(\LL, V)$. If $w(p)\ne 0$, one says that $p$
is an \emph{inflection point}. If $p$ is not an inflection point,
then $a_i(p)=i$ for all $i$, and vice versa. If we choose an
isomorphism $\psi\colon \LL \to \O_p$ and a generator $\pi \in
\m_p$, then the basis $s_0,\dots,s_n$ can be chosen in such a way
that $\psi (s_i) \equiv \pi^{a_i} \bmod {\m_p^{a_i+1}}$. We will
say that such a basis is \emph{adapted to $p$}. We will use the
following result (cf.\ \cite{EH, Proposition~1.1}):

\proclaim{Theorem 2.1 (\Plucker\ formula)}
If $C$, $\LL$, and $V$ are as above, then either all points of
$C$ are inflection points, or almost all points of $C$ are not
inflection points, and the following formula holds:
$$
  \sum_{p\in C} w(p)=(n+1)d+n(n+1)(g-1),
$$
where $d=\deg \LL$ and $g$ is the genus of $C$.
\endproclaim
If the characteristic is zero, then almost all points of $C$ are
not inflection points (cf.~\cite{GH}).

The \Plucker\ formula implies the following well-known
\proclaim{Proposition 2.1}
Suppose that no point of $C$ is an inflection point with respect
to the linear system $(C,\LL,V)$. Then this linear system defines
an isomorphism of $C$ onto a rational normal curve in $\PP^n$.
\endproclaim

\head 3. Abstract vanishing sequences
\endhead
For a non-negative integer $n$, denote by $\A_n$ the set of
strictly increasing sequences of $n+1$ non-negative integers:
$$
 \A_n=\set{(a_0,\dots,a_n)}
        0\le a_0< \dots < a_n
      \endset.
$$
Elements of $\A_n$ will be called \emph{abstract vanishing
sequences} of length~$n$.

For any integer $k\ge 1$ and any $A=(a_0,\dots,a_n)\in \A_n$,
put
$$
  v_k(A)=\left\{\,
         a_{i_1}+\dots+a_{i_k}\mid
         0\le i_1\le\dots\le i_k\le n
         \,\right\}.
$$
Sometimes we will assume that $v_0(a)=\{0\}$. By \emph{$k$-span}
of $A$ we will mean the cardinality of $v_k(A)$. We will denote
$k$-span of $A$ by $\Span_k(A)$.

The following easy observation plays the key role in the sequel.
\proclaim{Proposition 3.1}
If $(\LL, V)$ is a linear sytem on a smooth curve $C$, where
$\dim \LL=n+1$, and if $A=(a_0,\dots,a_n)$ is the vanishing
sequence at a point $p\in C$, then $\dim \im \left (\Sym^m V\to
H^0(\LL^{\otimes m}) \right) \ge \Span_m(A)$ for any $m\ge 1$.
\endproclaim
\demo{Proof}
Let $(s_0, \dots, s_n)$ be an adapted to $p$ basis of $V$. Among
$m$-fold products of $s_j$'s there are exactly $\Span_m(A)$ that
have different orders of vanishing at $p$. It is clear that all
these sections of $\LL^{\otimes m}$ are linearly independent,
whence the result.
\enddemo
\proclaim{Corollary 3.2}
If $X\subset \PP^n$ is an irreducible non-degenerate curve such
that $X= \phi(C)$, where $C$ is a smooth curve and $\phi$ is
defined by a linear system $(\LL, V)$ on $C$, and if $p\in C$
is a point with vanishing sequence $(a_0,\dots, a_n)$, then $\dim
\lspan{v_m(X)} \ge \Span_m(A) -1$.
\endproclaim

Now let us find out in what cases $m$-span of a sequence is
small.
\proclaim{Proposition 3.3}
Suppose that $A\in \A_n$. Then
\roster
\item"(i)"
$\Span_m(A)\ge mn+1$;
\item"(ii)"
$\Span_m(A)= mn+1$ if and only if $A$ is an arithmetic
progression.
\item"(iii)"
If $\Span_m(A)> mn+1$, then $\Span_m(A)\ge m(n+1)$.
\item"(iv)"
Suppose that $n\ge 3$. Then $\Span_m(A) = m(n+1)$ if and only if
either $a_n - a_{n-1} = 2(a_j- a_{j-1})$ for all $j\in [1; n-1]$
or $a_1 - a_0 =2(a_j-a_{j-1})$ for all $j\in [2;n]$.
\endroster
\endproclaim
\demo{Proof}
For $1\le i,j\le n$, put $b_{ij}=(m-i)a_0 +a_j +(i-1)a_n \in
v_m(A)$. It is clear that the following chain of inequalities
holds:
$$
\matrix
ma_0 & < & b_{11} < b_{12} <\dots < b_{1n} <{}\\
        && b_{21} < b_{22} <\dots < b_{2n} <{}\\
        && \hdotsfor 1\\
        && b_{m1} < b_{m2} <\dots < b_{mn}
\endmatrix
\tag 3.1
$$
Since this chain contains $mn+1$ elements of $v_m(A)$,
assertion (i) follows. Suppose now that $A$ is an arithmetic
progression with step $d$. Then the difference of any two
elements of $v_m(A)$ is a multiple of $d$, and any two
consecutive members of the chain~(3.1) differ by $d$.
Since $ma_0$ and $b_{mn}=ma_n$ are the least and the greatest
element of the set $v_m(A)$, we see that there are no elements in
$v_m(A)$ except for members of (3.1); this proves the
``if'' part of assertion~(ii).

To prove the ``if'' part of assertion~(iv), assume that $a_1-a_0
= 2d$ and $a_j - a_{j-1} = d$ for all $j >1$ (in the other case
the proof is analogous). Then the difference of any two elements
of $v_m(A)$ is again a multiple of $d$, and any two consecutive
members of the chain~(3.1) differ by $d$, except that
$$
   b_{21}-b_{1n}= b_{31}-b_{2n}=\dots = b_{m1}- b_{m-1,n}=2d.
$$
It is clear that
$$
  b_{jn}  < (m-j-1)a_0 + a_1 + a_{n-1} + (j-1)a_n < b_{j+1,1}
  \qquad\text{for all $j<m$.}
$$
If we insert the $m-1$ numbers $(m-j-1)a_0 + a_1 + a_{n-1}$ for
$j\in [1;m-1]$ into the chain~(3.1), we obtain a
strictly increasing sequence of length $m(n+1)$ such that the
difference of any pair of its consecutive elements equals $d$.
Since the difference of any two elements of $v_m(A)$ is divisible
by $d$, no extra element can be inserted into this sequence,
whence the cardinality of $v_m(A)$ equals $m(n+1)$.

To proceed, we need the following
\proclaim{Lemma 3.4}
\roster
\item"(i)"
If $(b_{k,n-1}; b_{k+1,1})\cap v_m(A)=\{b_{kn}\}$ for some $k\in
[1; m-1]$, then $a_n- a_{n-1} = a_1 - a_0$.
\item"(ii)"
If $a_{j+1}-a_j = a_1-a_0$ for some $j\in [2;n-1]$ and
$(b_{k,j-1}; b_{k,j+1}) \cap v_m(A) = \{b_{kj}\}$ for some $k\in
[1; m-1]$, then $a_j - a_{j-1} = a_1 - a_0$.
\endroster
\endproclaim
\demo{Proof of the lemma}
To prove (i), observe that the number $x= (m-k-1)a_0 + a_1 +
a_{n-1}+ (k-1)a_n \in v_m(A)$ belongs to the interval
$(b_{k,n-1}; b_{k+1,1})$; thus $x=b_{kn}$, whence $a_n - a_{n-1}
= a_1 - a_0$.

To prove (ii), observe that the hypothesis implies the equality
$$
  b_{k,j+1}=(m-k-1)a_0+a_1 + a_j+(k-1)a_n.
$$
Hence, the number $x= (m-k-1)a_0 +a_1 + a_{j-1} +(k-1)a_n$
belongs to the interval $(b_{k,j-1}; b_{k,j+1})$; thus $x=
b_{kj}$, whence $a_j - a_{j-1} = a_1 - a_0$.
\enddemo
To prove assertion~(iii) of Proposition~3.3 together with the
``only if'' part of~(ii), it suffices to prove that $A$ is an
arithmetic progression whenever $\Span_m(A) < m(n+1)$. Indeed, in
this case at most $m-2$ elements of $v_m(A)$ are not contained in
the chain~(3.1). Since for $k\in [1;m-1]$ the
intervals~$(b_{k,n-1}; b_{k+1,1})$ are disjoint, it follows that
$(b_{k,n-1}; b_{k+1,1}) \cap v_m(A) = \{ b_{kn} \}$ for some
$k\in [1;m-1]$, whence $a_n-a_{n-1}= a_1-a_0$ by virtue of
Lemma~3.4(i). The same argument shows that for each $i\in
[2;n-1]$ there exists a number $k_i \in [1; n]$ such that
$(b_{k_i,i-1}; b_{k_i,i+1}) \cap v_m(A) = \{ b_{k_i,i} \}$.
Applying Lemma~3.4(ii) $n-2$ times, we see that $a_{n-1}-a_{n-2}
= \dots = a_1-a_0$, i.e.\ $A$ is an arithmetic progression.

It remains to prove the ``only if'' part of~(iv).
First we do it for $m=2$:
\proclaim{Lemma 3.5}
Suppose that $A\in \A_n$, where $n\ge 3$, and $\Span_2(A)=2n+2$.
Then either $a_n - a_{n-1} = 2(a_j- a_{j-1})$ for all $j\in [1;
n-1]$ or $a_1 - a_0 =2(a_j-a_{j-1})$ for all $j\in [2;n]$.
\endproclaim
\demo{Proof of the lemma}
For any $A=(a_0, \dots,a_n) \in \A_n$, put $A'=(a_0,\dots,
a_{n-1})\in \A_{n-1}$ and $A''=(a_1,\dots, a_n)\in \A_{n-1}$.
First observe that
$$
\Span_2(A)\ge \Span_2(A')+2 \qquad\text{and }
\Span_2(A)\ge \Span_2(A'')+2.
\tag 3.2
$$
Indeed, since $A$ is a strictly increasing sequence, we see that
$a_{n-1}+a_n\ne 2a_n$, and none of these numbers equals
$a_i+a_j$, where $0\le i\le j <n$. Hence, $v_2(A)\setminus
v_2(A') \supset \{a_{n-1}+a_n, 2a_n\}$, whence the desired
inequality. The proof of the second inequality is analogous.

Now we proceed by induction on $n$. For $n=3$, the proof is
straightforward. Suppose now that the lemma is proved for all
$A\in \A_k$, where $k<n$, and that $n>3$. If $A\in \A_n$ and
$\Span_2(A)=2n+2$, then~(3.2) implies that $\Span_2(A')\le 2n$
and $\Span_2(A'')\le 2n$. Moreover, we see that one of these
inequalities is actually an equality (if this is not the case,
then part~(ii) of Proposition~3.3, which we have already proved,
implies that $A'$ and $A''$ are arithmetic progressions, whence
$A$ is an arithmetic progression and $\Span_2(A)=2n+1$, contrary
to the assumption). From now on we will assume that
$\Span_2(A')=2n$ (the argument for the case when
$\Span_2(A'')=2n$ is the same if one mirrors the indices with
respect to $n/2$). Observe that $\{a_{n-1}+a_n, 2a_n\} \cap
v_2(A')= \emptyset$; it follows from our assumptions that
$v_2(A)= v_2(A')\cup \{a_{n-1}+a_n, 2a_n\}$, whence
$a_{n-2}+a_n\in v_2(A')$. The only element of $v_2(A')$ that is
\emph{a priori} not less than $a_{n-2}+a_n$, is $2a_{n-1}$.
Hence, $a_{n-2}+a_n=2a_{n-1}$.

If we apply the induction hypothesis to $A'$, we see that either
$a_1,\dots,a_{n-1}$ is an arithmetic progression with step $d$
and $a_1-a_0=2d$, or $a_0,\dots,a_{n-2}$ is an arithmetic
progression with step $d$ and $a_{n-1} - a_{n-2}=2d$. In the
first case it is clear that $a_1,\dots,a_n$ is an arithmetic
progression with the same step $d$, and we are done. In the
second case it is easy to show that $\Span_2(A)=2n+3$, which
contradicts the hypothesis. This completes the proof of the
lemma.
\enddemo

To prove~(iv) in full generality, we may assume that $m>2$.
Denote by $\E$ the set of elements of $v_m(A)$ that do not belong
to the chain~(3.1). Suppose that $v_m(A) = m(n+1)$, i.e.\ $\#
\E=m-1$. I claim that each of the $m-1$ intervals~$(b_{k,n-1};
b_{k+1,1})$ contains exactly one element of $\E$. Indeed, assume
the converse; then Lemma~3.4(i) implies that $a_n-a_{n-1}=
a_1-a_0$. Moreover, arguing as in the proof of assertion~(iii)
above, we see that for each $i\in [2;n-1]$ there exists a number
$k_i \in [1; n]$ such that $(b_{k_i,i-1}; b_{k_i,i+1}) \cap
v_m(A) = \{ b_{k_i,i} \}$; applying Lemma~3.4(ii) $n-2$ times we
see that $A$ is an arithmetic progression--a contradiction.

Hence, for each $k\in [1;m-1]$ there exists a number $x_k\in
(b_{k,n-1}; b_{k+1,1})$, $x_j\ne b_{kn}$, and each of the
intervals $(b_{i-1,1};b_{i+1,1})$, $i\in [2;n-1]$ does not
contain elements of $v_m(A)$. Now Lemma~3.4(ii)
implies that $A$ is an arithmetic progression whenever
$a_n-a_{n-1}= a_1-a_0$. Hence, $a_n-a_{n-1}\ne a_1-a_0$, and for
any $k\in [1;n-1]$ the number $(m-k-1)a_0 + a_1 + a_{n-1}+
(k-1)a_n \in v_m(A)$ belongs to the interval $(b_{k,n-1};
b_{k+1,1})$ and does not coincide with $b_{kn}$. This shows that
$x_k=(m-k-1)a_0 + a_1 + a_{n-1}+ (k-1)a_n$. If is clear that
either $b_{k,n-1}<x_k <b_{kn}$ for all $k\in [1;n-1]$ (this is
the case if $a_1 -a_0 < a_n -a_{n-1}$), or $b_{kn} < x_k<
b_{k+1,1}$ (this is the case if $a_1 -a_0 > a_n -a_{n-1}$).

In the first case, there are exactly $2n+2$ elements of $v_m(A)$
in the interval $[b_{m-2,n};b_{mn}]$. Since a number $x\in
v_m(A)$ belongs to $[b_{m-2,n};b_{mn}]$ whenever $x=(m-2)a_n +
y$ with $y\in v_2(A)$, we see that $\Span_2(A)=2n+2$, and the
conclusion about $A$ follows from Lemma~3.5.

In the second case the argument is analogous if we consider
$v_m(A) \cap [ma_0; b_{1n}]$. This completes the proof.
\enddemo

\remark{Remark 3.6}
For a given $A\in \A_n$ and all $m\gg 0$, the number $\Span_m(A)$
is easily computable. An explicit formula is contained
in~\cite{B, Theorem 4.1} (see also Proposition~4.10 below).
\endremark

By way of an amusing application, we show how Proposition~3.3(iv)
yields a proof of the following well-known
\proclaim{Proposition 3.7}
If $C\subset \PP^n$ is an elliptic curve of degree $n+1$,
then $C$ has exactly $(n+1)^2$ distinct inflection
points, and the weight of each of these points equals~$1$ (or,
equivalently, the vanishing sequence is of the form
$(1,\dots,n-1,n+1)$).
\endproclaim
\demo{Proof}
It follows immediately from the \Plucker\ formula that
$\sum_{p\in C} w(p)=(n+1)^2$; the point is that there are
$(n+1)^2$ distinct inflection points. To prove this it suffices
to show that $w(p)=1$ for any inflection point $p\in C$. To that
end, denote by $A$ the vanishing sequence at $p$; since $C$ is
contained in $n(n-1)/2-1$ linearly independent quadrics, we see
that $\dim \im(H^0(\Op n(2)) \to H^0 (\O_C(2)))\le 2n+2$, whence
$\Span_2(A)\le 2n+2$ by virtue of Proposition~3.1. Since $a_0=0$
and $a_1=1$, we see, by virtue of Proposition~3.3(iv), that
$A=(0,1,\dots, n-1, n+1)$ and $w(p)=1$.
\enddemo
\remark{Remark 3.8}
Of course, this fact can be established directly by Riemann-Roch.
The $(n+1)^2$ inflection points are points of order $n+1$ on the
elliptic curve.
\endremark

\head 4. Adapted bases and filtrations
\endhead
If $R=k[X_0,\dots, X_n]$ is the ring of polynomials over a
field~$k$, then $R_m\subset R$ denotes the space of homogeneous
polynomials of degree~$m$. Let $A= (a_0,\dots, a_n)$ be an
abstract vanishing sequence. Suppose that $\xi =X_0^{k_0}\dots
X_n^{k_n} \in R$ is a monomial. Then the number $a_0k_0+a_1k_1+
\dots +a_n k_n$ is called \emph{weight} of $\xi$ and is denoted
by~$\wt^A(\xi)$ (or simply $\wt(\xi)$). Denote by $R^{[j]}
\subset R$ the $k$-subspace spanned by the monomials of
weight~$j$. Then $R^{[i]}R^{[j]} \subset R^{[i+j]}$ and
$R=\bigoplus R^{[j]}$. Thus, the $R^{[j]}$'s define a (\quasi)
grading of $R$; we will say that elements of $R^{[j]}$ are
\quasi\ of weight~$j$. If $R_m^{[j]}= R_m\cap R^{[j]}$, then
$R=\bigoplus R_m^{[j]}$; we will sometimes refer to elements of
$R_m^{[j]}$ as \emph{bihomogeneous} of bidegree $(m,j)$.

For $j,m\ge 0$ put
$$
  (J_m^A)^{[j]}=\set{f \in R_m^{[j]}}
  f(1,1,\dots,1)=0
  \endset.
$$
Put $J^A_m=\bigoplus_j (J_m^A)^{[j]}$, $(J^A)^{[j]}= \bigoplus_m
(J_m^A)^{[j]}$, and $J^A=\bigoplus J^A_m = \bigoplus
(J^A)^{[j]}$. Then $J^A \subset R$ is an ideal that is both
homogeneous and \quasi. Put $S^A= R/J^A$; the ring $S^A$ inherits
both homogeneous and \quasi\ grading from $R$. Homogeneous
(resp.\ \quasi) componenets of $S^A$ will be denoted by $S^A_m$
(resp.\ $(S^A)^{[j]}$). Dimensions of homogeneous components of
$S^A$ are easily computed (cf.\ \cite{B, Section 4}):
\proclaim{Proposition 4.1}
$\dim_k(S^A_m)= \Span_m(A)$.
\endproclaim
\demo{Proof}
It is clear from the definition that $\dim (J_m^A)^{[j]} =\max
(0,\dim R_m^{[j]} -1)$. Hence, $\dim S^A_m = \sum_j \dim
(S^A_m)^{[j]} = \sum_m (\dim R_m^{[j]} - \max(0, \dim R_m^{[j]}
-1)) = \# \set{j} \dim R_m^{[j]} \ne 0 \endset$. The latter
number equals $\Span_m(A)$, which completes the proof.
\enddemo
Some simple properties of rings $S^A$ are gathered in the
following
\proclaim{Proposition 4.2}
Let $A=(a_0,\dots,a_n)$ and $B=(b_0,\dots,b_n)$ be abstract
vanishing sequences.
\roster
\item"(i)"
If there exists a constant $c$ such that $b_i=a_i+c$ for all $i$,
then $S^A\cong S^B$.
\item"(ii)"
If there exists a constant $d$ such that $b_i=a_i\cdot d$ for all
$i$, then $S^A\cong S^B$.
\item"(iii)"
If $a_i=b_{n-i}$ for all $n$, then $S^A\cong S^B$.
\endroster
In all three cases we mean isomorphism of graded rings with
respect to \emph{homogeneous} grading.
\endproclaim
\demo{Proof}
Left to the reader.
\enddemo
\proclaim{Corollary 4.3}
Let $A=(a_0,\dots,a_n)$ and $B=(b_0,\dots,b_n)$ be abstract
vanishing sequences.
\roster
\item"(i)"
If there exists a constant $c$ such that $b_i=a_i+c$ for all $i$,
then $\Span_m(A)=\Span_m(B)$ for all $m$.
\item"(ii)"
If there exists a constant $d$ such that $b_i=a_i\cdot d$ for all
$i$, then $\Span_m(A)=\Span_m(B)$ for all $m$.
\item"(iii)"
If $a_i=b_{n-i}$ for all $n$, then $\Span_m(A)=\Span_m(B)$ for
all $m$.
\endroster
\endproclaim

Let $(\O, \m, k)$ be a discrete valuation ring. We will always
assume that $k$ is     embedded into $\O$. Choose a generator
$\pi \in \m$. Suppose that there exists an increasing sequence of
integers $a_0< a_1 <\dots <a_n$ and a sequence $s_0,s_1, \dots,
s_n \in \O$ such that $\ord(s_j)= a_j$ and $s_j\equiv \pi^{a_j}
\bmod {\m^{a_j+1}}$. Put $\tilde I=\set{f\in R} f(s_0,
\dots, s_n)=0\endset$. Denote by $I\subset \tilde I$ the
largest homogeneous ideal contained in $\tilde I$. Put $S=R/I$.
The \quasi\ grading induces a descending filtaration on $R$:
$$
F^i R=\sum_{j\ge i}R^{[j]}.
\tag 4.1
$$
Put $F^i I=F^iR \cap I$. The following proposition provides an
upper bound for dimensions of homogeneous components of $I$.
\proclaim{Proposition 4.4}
$\dim\left( F^j I_m/F^{j+1} I_m\right) \le \dim (J^A)_m^{[j]}$.
\endproclaim
\demo{Proof}
Let us construct an injection of $\left( F^j I_m/F^{j+1}
I_m\right)$ into $(J^A)_m^{[j]}$. To that end, associate to any
$f\in F^j I_m$ its quasihomogeneous component of weight~$j$. We
are to check two things:
\roster
\item"(i)"
The kernel of this morphism coincides with $F^{j+1}I_m$;
\item"(ii)"
The image of this morphism lies in $(J^A_m)^{[j]}$.
\endroster
Assertion~(i) is obvious; to prove~(ii), suppose that $f+g\in F^j
I_m$, where $f$ is bihomogeneous of weight $(m,j)$ and $g\in
F^{j+1} R$. Then $f(s_0,\dots,s_n)\equiv f(1,\dots,1)\pi^j
\bmod{\m^{j+1}}$. On the other hand, the equality $f(s_0,\dots,
s_n)+ g(s_0,\dots, s_n)=0$ implies that $f(s_0,\dots, s_n)\in
\m^{j+1}$. This proves~(ii) and the proposition.
\enddemo
\proclaim{Corollary 4.5}
$\dim S_m \le \dim S^A_m = \Span_m(A)$.
\endproclaim
This corollary suggests the following
\definition{Definition 4.1}
A system $(s_0,\dots, s_n)$ with vanishing sequence $A$ is called
\emph{$m$-maximal} if $\dim S_m = \dim S^A_m$.
\enddefinition
Here's how this definition works:
\proclaim{Proposition 4.6}
If a system $(s_0,\dots, s_n)$ is $m$-maximal and if $\sum_ {t\ge
m}J^A_t=RJ^A_m$, then the system $(s_0,\dots, s_n)$ is
$t$-maximal for all $t\ge m$.
\endproclaim
\demo{Proof}
Observe that the system $(s_0,\dots, s_n)$ is $m$-maximal if and
only if the homomorphisms $\left(F^j I_m/F^{j+1} I_m\right) \to
(J^A)_m^{[j]}$ are isomorphisms for all $j$. Denote by $q(f)\in
R_m$ the sum of \quasi\ components of the lowest degree for any
$f\in R_m$.

Consider now an element $\bar f\in J^A_t$, where $t > m$; we are
to prove that $\bar f=q(f)$, where $f\in I_t$. To that end,
observe that $\bar f=\sum g_i \bar h_i$, where $\bar h_i\in
J^A_m$, by hypothesis; on the other hand, $\bar h_i=q(h_i)$,
where $h_i\in I_m$, because the system $(s_0,\dots, s_n)$ is
$m$-maximal. It is clear that one can put $f=\sum g_i h_i$.
\enddemo
Another application of the notion of $m$-maximality is the
following
\proclaim{Proposition 4.7}
If a system $(s_0,\dots, s_n)$ is $m$-maximal and $\sum_{t\ge
m}J^A_t = RJ^A_m$, then $\sum_{t\ge m}I_t = RI_m$.
\endproclaim
\demo{Proof}
Suppose that $f\in F^j I_t\setminus F^{j+1} I_t$, where $t\ge m$.
I claim that $f=\varphi+ h$, where $\varphi \in R I_m$ and $h\in
F^{j+1} I_t$.

Indeed, it follows from the proof of Proposition~4.4 that
$q(f)\in (J^A_m)^{[j]}$, whence $q(f)=\sum g_i h_i$, where
$h_i\in J^A_m$. The observation at the beginning of the proof of
Proposition~4.6 shows that $h_i = q(u_i)$, where $u_i\in I_m$. If
$\varphi= \sum g_i u_i$, then it is clear that $\varphi\in R I_m$
and $h=f- \varphi \in F^{j+1} I_m$. This proves our claim.

The proposition follows by induction from what we have proved,
because \quasi\ weights of elements of $R_t$ are bounded.
\enddemo

Let us apply our results to linear systems on curves. Let $(\LL,
V)$ be a linear system on a smooth curve~$C$. Fix a point $p\in
C$; put $(\O, \m) = (\O_p, \m_p)$ and fix an isomorphism $\LL_p
\to \O$. If $(s_0, \dots, s_n)$ is an adapted basis of $V$, then
we can treat $s_j$'s as elements of $\O$ and apply the previous
results to the system $(s_0,\dots, s_n)$. Observe that $S$ is the
homogeneous coordinate ring of the curve $\phi(C)\subset \PP^n$,
where $\phi$ is the mapping defined by the linear system $(\LL,
V)$.

We will say that a linear system $(\LL, V)$ is \emph{$m$-maximal
at the point~$p$} whenever $(s_0, \dots, s_n)$ is $m$-maximal in
the sense of Definition~4.1.

Besides the curve $C$, we will consider a certain rational curve
$C^A$, which depends only on the vanishing sequence $A$. By
definition, $C^A= \Proj(S^A)\subset \PP^n =\Proj R$ (we mean
$\Proj$ with respect to \emph{homogeneous} grading).
\proclaim{Proposition 4.8}
If $k$ is algebraically closed, then $C^A =\overline {\Phi^A
(\Bbb A^1)} \subset \PP^n$, where $\Phi^A \colon \Bbb A^1
\to \PP^n$ acts by the formula $t\mapsto (t^{a_0}: t^{a_1}:\dots
:t^{a_n})$.
\endproclaim
\demo{Proof}
Obvious.
\enddemo
Degree an arithmetic genus of the rational curve~$C^A$ are easily
computable. Indeed, let $A\in \A_n$ be an abstract vanishing
sequence. Observe that, in view of Proposition~4.2,
one may assume that $a_0=0$ and g.c.d.\ of $a_j$'s equals $1$.
\proclaim{Proposition 4.9}
Under the above assumptions on $A$, degree of $C^A$ equals $a_n$
and arithmetic genus of $C^A$ equals $L_0 + L_\infty$, where
$L_0$ is the number of gaps of the semigroup $\left<a_1,\dots,
a_n\right>$ and $L_\infty$ is the number of gaps of the semigroup
$\left<a_n-a_{n-1}, \dots, a_n-a_1, a_n\right>$.
\endproclaim
\demo{Proof}
It follows from the hypothesis that there exist integers
$c_0,\dots, c_n$ such that $\sum c_i a_i=1$. Put $\psi(x_1,
\dots, x_n)=\prod x_i^{c_i}$. Then the rational map $\psi\colon
C^A\to \Bbb A^1$ is a birational inverse to $\Phi^A$. Hence,
$\deg C^A=a_n$ and $p_a(C^A) =\operatorname {length}
(\Phi_*\O_{\PP^1}/\O_{C^A})$, where $\Phi$ is the extension of
$\Phi^A$ to $\O_{\PP^1}$. Since $\psi$ is regular on $C^A$
outside $\{0,\infty\}$, this length is the sum of the lenghts of
stalks of $\Phi_*\O_{\PP^1} /\O_{C^A}$ at $\Phi(0)$ and $\Phi
(\infty)$. The stalk of $\Phi_*\O_{\PP^1} /\O_{C^A}$ at $\Phi(0)$
is isomorphic to $k[[t]]/k[[t^{a_1}, \dots, t^{a_n}]]$, and its
length clearly equals $L_0$. Ditto for $\Phi(\infty)$.
\enddemo
Now we can compute $\Span_m(A)$ for a given $A\in \A_n$ and $m\gg
0$. The following proposition is contained in~\cite{B, Theorem 4.1}.
\proclaim{Proposition 4.10}
Let $A\in \A_n$ be an abstract vanishing sequence. Suppose that
$a_0=0$ and g.c.d.\ of $a_j$'s equals $1$. Then for all $m\gg 0$
one has $\Span_m(A)=a_n\cdot m+ 1- L_0 - L_\infty$, where $L_0$
and $L_\infty$ are as in Proposition~4.9.
\endproclaim
\demo{Proof}
Proposition~4.1 implies that, for $m\gg 0$,
$\Span_m(A)=P(m)$, where $P$ is the Hilbert polynomial of the
curve $C^A$. Now Proposition~4.9 applies.
\enddemo

Proposition~4.6 yields the following
\proclaim{Proposition 4.11}
If $(\LL, V)$ is $m$-maximal for some $m$ and $\sum_{j\ge m}
J^A_m =R J^A_m$, then the degree and arithmetic genus of
$X=\phi(C)$ are equal to those of $C^A$.
\endproclaim
\demo{Proof}
Proposition~4.6 shows that $(\LL, V)$ is $t$-maximal for all
$t\ge m$. Hence, Hilbert polynomials of $X$ and $C^A$ coincide.
\enddemo

\head 5. On equations defining certain monomial curves
\endhead
To apply Proposition~4.11, one should know the degrees of
generators of the homogeneous ideal $J^A$, or, equivalently, of
the homogeneous ideal of the curve $C^A$. In general, it is not
clear how these degrees depend on $A$. The following proposition
shows that in the simplest cases $J^A$ is generated by quadrics.
\proclaim{Proposition 5.1}
Let $A\in \A_n$ be an abstract vanishing sequence. Suppose that
one of the following conditions holds:
\roster
\item"(i)"
$n\ge 2$ and $A$ is an arithmetic progression;
\item"(ii)"
$n\ge 3$ and either $a_n - a_{n-1} = 2(a_j- a_{j-1})$ for all
$j\in [1; n-1]$ or $a_1 - a_0 =2(a_j-a_{j-1})$ for all $j\in
[2;n]$.
\endroster
Then $J^A=R J^A_2$.
\endproclaim
\proclaim{Corollary 5.2}
If $A\in \A_n$ is as in Proposition~5.1, then
$\sum_{j\ge m} J^A_m =R J^A_m$ for all $m\ge 2$.
\endproclaim
\remark{Remark 5.3}
If $n=2$ in the case~(ii), then $C^A$ is the plane cuspidal cubic
and is not defined by quadratic equations.
\endremark
\remark{Remark 5.4}
The curve $C^A\subset \PP^n$ is a rational normal curve in the
case~(i) and a linearly normal curve of arithmetic genus~$1$ in
the case~(ii). It is well known that the homogeneous ideals of
such curves are generated by its components of degree~$2$ (this
is especially true for the case~(i)). However, we give a proof
since no suitable reference for the case~(ii) is known to the
author.
\endremark

The rest of this section is devoted to the proof of
Proposition~5.1. We start with some general
observations. Consider an arbitrary $A=(a_0,\dots,a_n)$. It is
clear from the definitions that $J^A_m \subset k[X_0, \dots,
X_n]$ is generated, as a linear space, by differences $\xi-\eta$,
where $\xi$ and $\eta$ are monomials of degree $m$ and equal
weight (recall that weight of a monomial~$\xi=X_0^{d_0} \dots
X_n^{d_n}$ is the number $\wt^A(\xi)=\sum a_id_i$). For an
integer $t\ge 2$, we will say that monomials $\xi= X_0^{d_0}
\dots X_n^{d_n}$ and $\eta= X_0^{e_0} \dots X_n^{e_n}$ are
\emph{$t$-neighbors} whenever $\deg\xi= \deg \eta$,
$\wt^A(\xi)=\wt^A(\eta)$, and $\# \set{i\in [0;n]} d_i\ne
e_i \endset=t$. We will say that the monomials $\xi$ and $\xi'$
are \emph{$t$-equivalent} whenever there exists a chain of
monomials $\xi=\xi_0, \xi_1,\dots, \xi_l=\xi'$ such that for all
$j\in [1;l]$ the monomials $\xi_j$ and $\xi_{j-1}$ are
$t$-neighbors. The following propositions follow immediately
from the above definitions.
\proclaim{Proposition 5.5}
$J^A_m = R_{m-t}J^A_t$ if and only if any two monomials $\xi$ and
$\eta$ of degree~$m$ and equal weight are $t$-equivalent.
\endproclaim
\proclaim{Proposition 5.6}
If the monomials $\xi$ and $\eta$ are $t$-equivalent, then
$\lambda\xi$ is $t$-equivalent to $\lambda\eta$ for any monomial
$\lambda$.
\endproclaim
Suppose that $\xi= X_0^{d_0} \dots X_n^{d_n}$ is a monomial; let
us say that the \emph{support of $\xi$} is the set $\supp(\xi)
=\set{i\in [0;n]} d_i\ne 0\endset$. We will say that two
subsets $U,V$ of a segment $[k;l]$ are \emph{interlaced} iff
either there exist numbers $p',p'' \in U$ and $q\in V$ such that
$p'< q <p''$ or there exist numbers $p\in U$ and $q',q'' \in V$
such that $q' < p <q''$.
\proclaim{Lemma 5.7}
If $\xi$ and $\eta$ are monomials of equal weight, then
$\supp(\xi)$ and $\supp(\eta)$ are interlaced.
\endproclaim
\demo{Proof}
Assume the converse. If $\xi=X_0^{d_0} \dots X_n^{d_n}$ and
$\eta= X_0^{e_0} \dots X_n^{e_n}$, then we may suppose without
loss of generality that $d_i < e_i$ for all~$i$. Hence, $\sum a_i
d_i< \sum a_i e_i$, which contradicts the hypothesis.
\enddemo

To proceed with the proofs, it will be convenient to consider a
certain ``mathematical game.'' This game is played on a board
that is
a rank of $n+1$ squares numbered by integers $0,1,\dots, n$ from
left to right. Any number of pieces can stand on each square.
If we assign a position with $d_i$ pieces on the $i$-th square
to each monomial $\xi= X_0^{d_0} \dots X_n^{d_n}$, we get a
\hbox{1-1} correspondence between the set of monomials in
$X_0,\dots, X_n$ of degree $m$ and the set of positions in this
game with $m$ pieces on the board.

The pieces can be moved according to the following rule: If
$a_i+a_j= a_{i'} + a_{j'}$ and there are pieces on the $i$-th
and $j$-th squares, then both these pieces can be moved
simultaneously to the $i'$-th and $j'$-th square respectively.
This rule does not exclude the cases $i=j$ or $i'=j'$ (of
course, if $i=j$, then there should be at least two pieces on
the $i$-th square).

Using this language, one can restate the definition of
$2$-equivalence as follows:
\proclaim{Proposition 5.8}
Two monomials of equal degree are $2$-equivalent if and only if
the corresponding positions in the game can be joined by a
sequence of moves.
\endproclaim

Now we can resume the proofs.
\proclaim{Lemma 5.9}
Suppose that there exist integers $k$ and $l$ such that $0\le k <
l\le n$ and the sequence $(a_k, a_{k+1}, \dots, a_l)$ is an
arithmetic progression. If $\xi$ and $\eta$ are two monomials
such that the sets $\supp(\xi)\cap [k;l]$ and $\supp(\eta)\cap
[k;l]$ are interlaced, then there exist monomials $\xi'$ and
$\eta'$ such that $\xi'$ is $2$-equivalent to $\xi$, $\eta'$ is
$2$-equivalent to $\eta$, and $\supp(\xi')\cap \supp(\eta') \ne
\emptyset$.
\endproclaim
\demo{Proof}
Without loss of generality we may assume that there exist
integers $p',p'' \in \supp(\xi)$ and $q\in \supp(\eta)$ such that
$k\le p' < q < p''\le l$. Consider a position in our game that
corresponds to the monomial $\xi$. There are pieces on the
squares~$p'$ and $p''$; since the restriction of the sequence
$(a_0,\dots, a_n)$ to the segment $[p'; p'']$ is an arithmetic
progression, we may start moving both these pieces towards each
other, to one square each at a time, until they are on the same
or adjacent square(s). En route, one of the pieces will step on
the $q$-th square; if $\xi'$ is the monomial corresponding to the
position at that moment, then $\xi'$ is $2$-equivalent to $\xi$
and $\supp(\xi') \cap \supp (\eta) \ne \emptyset$.
\enddemo

\demo{Proof of Proposition~5.1(i)}
In view of Proposition~5.5 it suffices to prove that any two
monomials of degree $m$ and equal weight are $2$-equivalent. Let
us do it by induction on $m$.

If $m=2$, there is nothing to prove. Suppose that our assertion
is proved for all $m' < m$, and let $\xi$ and $\eta$ be two
monomials of degree $m$ and equal weight. Lemma~5.7 implies that
$\supp(\xi)$ and $\supp(\eta)$ are interlaced; since $A$ is an
arithmetic progression, Lemma~5.9 implies that there exist
monomials $\xi'$ and $\eta'$ such that $\xi'$ is $2$-equivalent
to $\xi$, $\eta'$ is $2$-equivalent to $\eta$, and $\supp(\xi')
\cap \supp(\eta') \ne \emptyset$. Hemce, there exist monomials
$\bar \xi$, $\bar \eta$ and $\lambda$ such that $\xi' =\lambda
\bar \xi$, $\eta' =\lambda \bar \eta$, and $\deg \lambda >0$.
Since weights of $\bar \xi$ and $\bar \eta$ are equal, the
induction hypothesis implies that $\bar \xi$ and $\bar \eta$ are
$2$-equivalent, whence $\xi$ and $\eta$ are $2$-equivalent by
virtue of Proposition~5.6. This completes the proof.
\enddemo

\demo{Proof of Proposition~5.1(ii)}
Proposition~4.2 shows that we may assume that
$A=(0,1, \dots, n-1, n+1)$, which we will do further on.
Proceeding by induction as in the previous proof, one can see
that it suffices to prove the following claim:
\block
If $\xi$ and $\eta$ are two monomials of degree $m>2$ and equal
weight, then there exist monomials $\xi'$ and $\eta'$ such that
$\xi'$ is $2$-equivalent to $\xi$, $\eta'$ is $2$-equivalent to
$\eta$, and $\supp(\xi')\cap \supp(\eta')\ne \emptyset$.
\endblock
To prove this claim, observe that $\supp(\xi)$ and $\supp(\eta)$
are interlaced by virtue of Lemma~5.7. If $n \notin
\supp(\xi) \cup \supp(\eta)$, then we are done by
Lemma~5.9. Hence, we may assume that $n\in \supp(\xi)$
and $\supp(\xi)\cap \supp(\eta)= \emptyset$ (in particular,
$\supp(\eta) \not \ni n$).

Now the game comes to help. Place $m$ white pieces and $m$
black pieces on the board in such a way that the position of
white (resp.\ black) pieces corresponds to the monomial $\xi$
(resp.\ $\eta$). We are to show that one can move white and\slash
or black pieces according to the rules of the game in such a way
that one of the pieces steps on a square occupied by a piece of
the opposite color.

If the sets $\supp(\xi)\cap [0;n-1]$ and $\supp(\eta)\cap
[0;n-1]$ are interlaced, then we can make white and black pieces
meet using the method from the proof of Lemma~5.9.
Hence, in the sequel we may and will assume that $\supp(\xi)\cap
[0;n-1]$ and $\supp(\eta)\cap [0;n-1]$ are \emph{not} interlaced.
It follows from this assumption that for any $i\in \supp(\xi)
\setminus \{n\}$ and $j\in \supp(\eta)$ we have $i<j$.

If $\supp(\xi) \supset \{n,j\}$, where $j\le n-3$, we can make
the following move: a white piece jumps from $j$ to $j+2$, and
another white piece jumps from $n$ to $n-1$. Let us repeat such a
move (with various $j$'s maybe) while it is possible. When this
process is finished, one of the following assertions holds:
\roster
\item"1."
There are no white pieces left on the $n$-th square.
\item"2."
There are no white pieces left on the squares whose number
is less than $n-2$, and there is a white piece on the $n$-th
square.
\endroster
Denote the monomial corresponding to the new position of white
pieces by $\xi'$. In the first case, we can finish off the proof
by applying Lemma~5.9 to $\xi'$ and $\eta$. In the second case,
all white pieces are concentrated on the $(n-2)$-th and $n$-th
squares. If $\supp(\xi') \cap [0;n-1]$ is interlaced with
$\supp(\eta)$, we can again make white and black piece meet by
the same method as in the proof of Lemma~5.9. Hence, we may
assume that these two sets are not interlaced and disjoint. This
is possible only if all black pieces are concentrated on the
$(n-1)$-th square. Since $n\ge 3$ by the hypothesis, we can make
the following move: one black piece jumps from $n-1$ to $n$, and
another jumps from $n-1$ to $n-3$. After this move, black and
white pieces  will meet on the $n$-th square. This completes the
proof.
\enddemo

\head 6. Proofs of main results
\endhead
\demo{Proof of Theorem~1.4}
If $X\subset \PP^n$ is an irreducible non-degenerate curve, then
$X$ is the image of a birational morphism $f\colon C\to \PP^n$
defined by a linear system $(\LL, V)$, where $C$ is a smooth
curve. If $p\in C$ is any point with the vanishing sequence
$A=(a_0,\dots,a_n)$, then Proposition~3.3(i) shows that
$\Span_m(A)\ge mn+1$, whence $\dim \im (\Sym^m (V)\to
H^0(\LL^{\otimes m})) \ge mn+1$ by Proposition~3.1. The first
part of the theorem follows immediately from this inequality.

The fact that the bound is attained for rational normal curves is
trivial. Hence, we can assume for the sequel that $h^0(\idsheaf
X(m)) = \binom {m+n}n-mn-1$, or, equivalently, that $\dim \im
(\Sym^m (V)\to H^0(\LL^{\otimes m})) = mn+1$, and we are to prove
that $X$ is a rational normal curve.

To that end, pick a point $p\in C$ and denote by $A=(a_0,\dots,
a_n)$ its vanishing sequence at $p$. It follows from
Proposition~3.1 that $\Span_m(A)\le mn+1$. Now
Proposition~3.3(i,ii) shows that $A$ is an arithmetic
progression and that $(\LL, V)$ is $m$-maximal at $p$.

Hence, $C^A$ is the normal rational curve of degree $n$ in
$\PP^n$. Proposition~4.11 and Corollary~5.2
show that $X$ has degree $n$ and arithmetic genus $0$ as well,
whence $X$ is rational normal curve.
\enddemo

\demo{Proof of Theorem~1.5}
Consider the linear system~$(\LL, V)$, where $\LL=\O_X(1)$ and
$V=\im (H^0(\O_{\PP^n}(1) \to H^0(\O_X(1))$. If this linear
system has no inflection points, then $X$ is a normal rational
curve by Proposition~2.1, whence $h^0(\idsheaf X(m))=
\binom {m+n}n-mn-1$, which contradicts the hypothesis. Hence, we
may assume that there is an inflection point $p\in X$. Denote its
vanishing sequence by $A=(a_0,\dots, a_n)$. It is clear that
$a_0=0$ and $a_1=1$, whence $A$ is not an arithmetic progression.

Proposition~3.3(iii) implies that $\Span_m(A)\ge
m(n+1)$, and Proposition~3.1 implies that $\dim \im
\left (\Sym^m V\to H^0(\LL^{\otimes m}) \right) \ge m(n+1)$,
whence the first assertion of the theorem.

To prove the second assertion, suppose that $h^0(\idsheaf X(m))
= \binom{m+n}n - m(n+1)$. Propositions~3.3(iv)
and~3.1 imply that $A$ is of the
form~$(0,1,\dots,n-1,n+1)$ and the linear system $(\LL, V)$ is
$2$-maximal at $p$. Corollary~5.2 shows that
Proposition~4.11 applies, whence the degree and
arithmetic genus of $X$ equal those of $C^A$. Since $C^A$ is a
cuspidal rational curve of degree $n+1$ and arithmetic genus~$1$,
we are done.
\enddemo

\demo{Proof of Theorems~1.2 and 1.3}
We will proceed by induction on $\dim X$. If $\dim X=1$, then our
theorem are just the $m=2$ case of Theorems~1.4
and~1.5. If $\dim X>1$, consider a generic hyperplane
section $Y\subset X$. There is an exact sequence
$$
  0\to \idsheaf X(1)\to \idsheaf X(2)\to \idsheaf Y(2)\to 0.
$$
Since $X$ is non-degenerate, this sequence yields the inequality
$h^0 \idsheaf X(2)\le h^0 \idsheaf Y(2)$, and the theorems follow
from the induction hypothesis.
\enddemo

\enddocument